\documentclass[draft,11pt,a4paper]{article}
\usepackage{amssymb,amsmath,theorem}
\usepackage[includeheadfoot,margin=1in]{geometry}

\allowdisplaybreaks[2]

\newcommand\eq{\leftrightarrow}
\newcommand\nequiv{\not\equiv}
\newcommand\LOR{\bigvee}
\newcommand\ET{\bigwedge}
\newcommand\model{\vDash}
\newcommand\nmodel{\nvDash}
\newcommand\lang{\mathcal L}
\newcommand\leqg{\lesseqgtr}

\newcommand\fii{\varphi}
\newcommand\ep{\varepsilon}

\newcommand\p[1]{\langle#1\rangle}
\newcommand\lh[1]{\lvert#1\rvert}
\newcommand\abs[1]{\lh{#1}}
\newcommand\bez{\smallsetminus}
\newcommand\sset{\subseteq}
\newcommand\ssset{\subsetneq}
\newcommand\sSset{\supsetneq}
\newcommand\nul{\varnothing}
\newcommand\res{\mathbin\restriction}
\newcommand\up{\uparrow}

\newcommand\fl[1]{\lfloor#1\rfloor}
\newcommand\delim[4]{\ifx X#3X\left#1#4\right#2\else\csname#3l\endcsname#1#4\csname#3r\endcsname#2\fi}
\newcommand\Fl[2][]{\delim\lfloor\rfloor{#1}{#2}}
\newcommand\cl[1]{\lceil#1\rceil}
\newcommand\CL[2][]{\delim\lceil\rceil{#1}{#2}}
\newcommand\tsum{\textstyle\sum}

\DeclareMathOperator\bit{bit}
\DeclareMathOperator\card{card}
\DeclareMathOperator\Th{Th}
\DeclareMathOperator\lcm{lcm}

\newcommand\cxt[1]{\mathrm{#1}}
\newcommand\tc{\cxt{TC}^0}

\newcommand\thry[1]{\mathsf{#1}}
\newcommand\delt{\Delta^b_}
\newcommand\delz{\Delta_0}
\newcommand\Sig{\Sigma^B_}
\newcommand\io{\thry{IOpen}}
\newcommand\idz{\thry{I\delz}}
\newcommand\texp{\thry{EXP}}
\newcommand\idze{\idz+\texp}
\newcommand\idzw{\idz+\thry{\Omega_1}}
\newcommand\dicr{\thry{\delt1\text-CR}}
\newcommand\vtc{\thry{VTC^0}}
\newcommand\GGG{\thry{3G}}
\newcommand\GG{\thry{2G}}

\newcommand\N{\mathbb N}
\newcommand\Q{\mathbb Q}
\newcommand\Z{\mathbb Z}
\newcommand\RR{\mathbb R}
\newcommand\strc[1]{\mathbf{#1}}
\newcommand\MN{\strc N}
\newcommand\ML{\strc L}
\newcommand\MZ{\strc Z}
\newcommand\MQ{\strc Q}
\newcommand\MR{\strc R}
\newcommand\MZL{\MZ_\ML}
\newcommand\MQL{\MQ_\ML}
\newcommand\MRL{\MR_\ML}
\newcommand\sM{\mathfrak M}

\newcommand\ob[1]{\overline{#1}}

\newcommand\noproof{\leavevmode\unskip\hskip.75em\vadjust{}\nobreak\hfill$\Box$\par}
\newenvironment{Pf}
  {\par\noindent\textit{Proof:}\hskip.75em\ignorespaces}
  {\noproof\pagebreak[2]\vskip\medskipamount\ignorespacesafterend}

\theoremstyle{plain}
\newtheorem{Thm}{Theorem}[section]
\newtheorem{Cor}[Thm]{Corollary}
\newtheorem{Lem}[Thm]{Lemma}
\newtheorem{Que}[Thm]{Question}
\theorembodyfont\upshape
\newtheorem{Def}[Thm]{Definition}
\newtheorem{Rem}[Thm]{Remark}
\newtheorem{Exm}[Thm]{Example}

\author{Emil Je\v r\'abek\\[\medskipamount]
Institute of Mathematics, Czech Academy of Sciences\\
\small \v Zitn\'a 25,
115\:67 Praha 1,
Czech Republic,
email: \texttt{jerabek@math.cas.cz}}

\title{Models of $\vtc$ as exponential integer parts}

\begin{document}
\maketitle

\begin{abstract}
We prove that (additive) ordered group reducts of nonstandard models of the bounded arithmetical theory $\vtc$ are
recursively saturated in a rich language with predicates expressing the integers, rationals, and logarithmically
bounded numbers. Combined with our previous results on the construction of the real exponential function on completions
of models of~$\vtc$, we show that every countable model of~$\vtc$ is an exponential integer part of a real-closed
exponential field.


\end{abstract}

\section{Introduction}\label{sec:introduction}

A classical result of Shepherdson~\cite{sheph} characterizes models of the arithmetical theory $\io$ as integer parts
(IP) of real-closed fields. Conversely, every real-closed field has an integer part by Mourgues and
Ressayre~\cite{mou-res}. Ressayre~\cite{ress:eip} introduced an analogous notion for exponential fields, namely
exponential integer parts of real-closed exponential fields (RCEF; the definition includes the growth axiom
$\exp(x)>x$). He proved that every RCEF has such an exponential IP. Here we are interested in the converse problem:
\begin{Que}\label{que:eip}
What ordered rings are exponential IP of real-closed exponential fields?
\end{Que}
It is not surprising that every model of $\idze$ is an exponential IP of a RCEF, though surprisingly, this does not
seem to be well established in the literature: the results we are aware of are that every model of $\idze$ is an
exponential IP of a real-closed field admitting so-called left exponential by Carl, D'Aquino, and Kuhlmann~\cite{cdk},
and every model of $\thry{PA}$ is an exponential IP of a RCEF by Krapp~\cite{krapp:phd}; see also Carl and
Krapp~\cite{carl-krapp}. In any case, we prove that every model of $\idze$ is an exponential IP of a RCEF below
(Corollary~\ref{cor:idz+exp}).

However, our main interest in this paper are models of weak theories of arithmetic where integer exponentiation
is not total. The definition of exponential IP does not require the field exponential to extend the usual
integer exponential function as considered in theories of arithmetic, yet we might wonder whether its growth and
algebraic properties perhaps force the totality of integer exponentiation nonetheless, or at least, whether they imply
some nontrivial first-order consequences of $\idze$. We essentially give a negative answer to both questions: the
first-order consequences of being an exponential IP of a RCEF are contained in~$\vtc$---a weak subtheory of bounded
arithmetic---and more specifically, every countable model of~$\vtc$ is an exponential IP of a RCEF (Theorem~\ref{thm:rcef}).

Uniform $\tc$ is a small computational complexity class below logarithmic space and polynomial time. It can be thought
of as the complexity of basic arithmetic operations: integer (and rational) $+$, $-$, $\cdot$, $/$, and $<$ are
$\tc$-computable, with $\cdot$ and $/$ being $\tc$-complete under a suitable notion of reduction. $\vtc$ is the basic
theory of $\tc$-computable functions, similar to how $\thry{PRA}$ is the basic theory of primitive recursive functions.
$\vtc$ as introduced by Nguyen and Cook \cite{ngu-cook} is a two-sorted theory of bounded arithmetic in the style of
Zambella~\cite{zamb:notes}, but for the purposes of this paper it can be identified with the one-sorted theory
$\dicr$ of Johannsen and Pollett~\cite{joh-pol:d1cr}, which is a weak fragment of
$\idzw$. Any model $\sM$ of $\vtc$ or~$\dicr$ induces a discretely ordered ring~$\MZ^\sM$ (the ``integers'' of~$\sM$),
the fraction field~$\MQ^\sM$ of~$\MZ^\sM$ (the ``rationals'' of~$\sM$), and the completion $\MR^\sM$ of~$\MQ^\sM$ (the
``reals'' of~$\sM$).

Most of the hard work was done in Je\v r\'abek~\cite{ej:vtcanal}, where it is shown that for any model $\sM\model\vtc$,
the field of reals~$\MR^\sM$ (which is real-closed by \cite{ej:vtc0iopen,ej:vtcimul}) carries a well-behaved
analogue of the real exponential function. However, if $\sM\nmodel\texp$, this exponential is not total: it
is only defined on $\MRL^\sM$, the logarithmically bounded reals. We overcome this problem by proving,
for countable~$\sM$, that $\p{\MR^\sM,\MZ^\sM,+,{<}}$ and $\p{\MRL^\sM,\MZL^\sM,+,{<}}$ are isomorphic, and one can
choose the isomorphism such that the resulting exponential on~$\MR^\sM$ satisfies the growth axiom $\exp(x)>x$.

The main part of our argument is to show that for any nonstandard $\sM\model\vtc$, the structure
$\p{\MQ^\sM,\MZ^\sM,\MQL^\sM,+,{<}}$ is recursively saturated (Theorem~\ref{thm:vtc-rec-sat}), which is a result of
independent interest. This is a continuation of a line of research showing that tame structures interpretable in
nonstandard models of sufficiently strong arithmetic are recursively saturated: in particular, the additive reduct
(corresponding to our $\p{\MZ^\sM,+,{<}}$) of a nonstandard $\sM\model\thry{IE_1}$ is recursively saturated by
Wilmers~\cite{wilmers}, following up on~\cite{lessan:phd,jen-ehr:prob,cmw}; for structures of another kind, if $\sM$ is
a nonstandard model of a suitable arithmetic, then any real-closed field with IP $\sM$ (such as our
$\p{\MR^\sM,+,\cdot}$) is recursively saturated by \cite{dks,ej-kol:realclosures}.

The paper is organized as follows. After this Introduction, some preliminary definitions and notation are summarized in
Section~\ref{sec:preliminaries}. Section~\ref{sec:real-expon-models} reviews what follows from the results
of~\cite{ej:vtcanal} and what is missing. In Section~\ref{sec:theory-three-groups}, we axiomatize the theory of
$\p{\MQ^\sM,\MZ^\sM,\MQL^\sM,+,{<}}$ (the ``theory of three groups''~$\GGG$, Definition~\ref{def:3g}), and prove a quantifier
elimination result for this theory (Theorem~\ref{thm:qe-3g}). In Section~\ref{sec:recurs-satur-models}, we characterize
recursive saturation of models $\p{Q,Z,L,+,{<}}\model\GGG$ in terms of recursive saturation of the $\p{Q,Z,+,{<}}$
reducts (Theorem~\ref{thm:3g-rec-sat}). We prove our main results (Theorems \ref{thm:vtc-rec-sat} and~\ref{thm:rcef}) in
Section~\ref{sec:ggg-reducts-models}, and conclude the paper with some open problems in Section~\ref{sec:conclusion}.

\section{Preliminaries}\label{sec:preliminaries}

In this paper, all groups are assumed to be abelian, and all ordered algebraic structures are assumed to be totally
ordered. In particular, an \emph{ordered group} is a structure $\p{G,+,0,{<}}$ such that $\p{G,+,0}$ is an
abelian group, and $<$ is a total order on~$G$ such that $x\le y\implies x+z\le y+z$ for all $x,y,z\in G$. We denote
the set of positive elements of $G$ as~$G_{>0}$. A subset $X\sset G$ is \emph{convex} if $[x,y]\sset X$ for
all $x,y\in X$ such that $x\le y$, where $[x,y]$ denotes the closed interval $\{z\in G:x\le z\le y\}$. If $x=y$, the
interval $[x,y]$ is \emph{degenerate}. We define also open and half-open intervals $(x,y)$, $[x,y)$, $(x,y]$ as usual,
including unbounded intervals with endpoints in $G\cup\{-\infty,+\infty\}$ where $-\infty<x<+\infty$ for all $x\in G$.

Any nontrivial ordered group is either dense(ly ordered) or \emph{discrete}, meaning that $G_{>0}$ has a least element
(often denoted $1$). An \emph{integer part} (\emph{IP}) of an ordered group~$G$ is a discrete subgroup $Z\sset G$ with a least
positive element~$1$ such that every $x\in G$ is within distance~$1$ from an element $n\in Z$ (i.e., $\abs{x-n}\le1$,
where $\abs x=\max\{x,-x\}$). Then for every $x\in G$, there is a unique element $\fl x\in Z$ such that $\fl x\le x<\fl
x+1$; we also write $\{x\}=x-\fl x$.

An \emph{ordered ring} is a structure $\p{R,+,\cdot,0,1,{<}}$ such that $\p{R,+,\cdot,0,1}$ is a commutative ring,
$\p{R,+,0,{<}}$ is an ordered group, and $x\le y\implies xz\le yz$ for all $x,y\in R$ and $z\in R_{>0}$. An ordered
ring is \emph{discrete} if $1=\min R_{>0}$. An \emph{integer part} (\emph{IP}) of an ordered ring $R$ is a discrete subring
$Z\sset R$ which is an IP of its additive group. An \emph{ordered field} is an ordered ring that is a field. A
\emph{real-closed field} is an ordered field $R$ with no proper algebraic ordered field extension; equivalently,
every $f\in R[x]$ of odd degree has a root in~$R$, and every $a\in R_{>0}$ has a square root in~$R$; also equivalently,
$R$ is elementarily equivalent to $\p{\RR,+,\cdot,0,1,{<}}$

Shepherdson~\cite{sheph} proved that $R$ is an IP of a real-closed field iff $R\model\io$, where $\io$ is the
theory of discrete ordered rings augmented with the induction schema
\[\fii(0,\vec y)\land\forall x\:\bigl(\fii(x,\vec y)\to\fii(x+1,\vec y)\bigr)\to\forall x\ge0\:\fii(x,\vec y)\]
for open (= quantifier-free) formulas~$\fii$.

An \emph{(ordered) exponential field} is an ordered field $R$ endowed with an ordered group isomorphism
$\exp\colon\p{R,+,0,{<}}\to\p{R_{>0},\cdot,1,{<}}$. Following Ressayre~\cite{ress:eip}, a \emph{real-closed exponential
field} is an exponential field $\p{R,+,\cdot,0,1,{<},\exp}$ which is real-closed and satisfies
$\exp(1)=2$ and%
\footnote{Ressayre actually demands ``$\exp(x)>x^n$ for all $x$ somewhat larger than $n$'', where $n$ presumably
refers to standard natural numbers. This follows from our formulation, since $\exp(x)=\exp(x/2n)^{2n}>(x/2n)^{2n}\ge
x^n$ as long as $x\ge(2n)^2$ (this can be improved). On the other
hand, it is easy to see that if $\exp(x)>x$ holds for all $x\ge m\in\N$, then it holds for all $x\in R$, thus our axiom
is equivalent to Ressayre's formulation.}
$\exp(x)>x$ for all $x\in R$. An \emph{exponential integer part} of an exponential field
$\p{R,\exp}$ is an IP $Z\sset R$ such that $Z_{>0}$ is closed under $\exp$.
Ressayre shows that every real-closed exponential field has an exponential IP (this is further elaborated
in~\cite{dkkl:rcef}).

Every ordered field $F$ has a \emph{completion $\hat F$} that can be described in several equivalent
ways. One way using only
the basic structure of ordered fields is as follows (cf.~\cite{scott-cof}). A \emph{cut} in $F$ is a pair $\p{A,B}$ of
sets such that $F=A\cup B$, $\inf\{b-a:b\in B,a\in A\}=0$, and $A$ has no largest element; $F$ is \emph{complete} if
$\min B$ exists for every cut $\p{A,B}$. The \emph{completion} of $F$ is a complete ordered field
$\hat F$ such that $F$ is a dense subfield of $\hat F$ (i.e., every non-degenerate interval of $\hat F$
intersects $F$). The completion of $F$ is unique up to $F$-isomorphism; it can be explicitly constructed by endowing
the set of all cuts of~$F$ with suitable structure.

We will most often use a topological description of~$\hat F$ (see~\cite{wer:topf}). The interval topology makes $F$ a
topological field, and therefore a uniform space\footnote{We require all uniform spaces and topological groups to be
Hausdorff.} with a fundamental system of entourages $\mathcal U=\{U_\ep:\ep\in F_{>0}\}$, where
$U_\ep=\{\p{x,y}\in F^2:\abs{x-y}\le\ep\}$. $F$ is \emph{complete} as a uniform space if every Cauchy net in $F$
converges.
The \emph{completion} of $F$ is a complete uniform space $\hat F$ such that $F$ is a (topologically) dense subspace
of~$\hat F$; it is again unique up to $F$-isomorphism. The key property of $\hat F$ is that every uniformly continuous
function from $F$ to a complete uniform space~$S$ extends uniquely to a uniformly continuous function $\hat F\to S$.
The ring operations on~$F$ extend to continuous operations on~$\hat F$ that make it a topological ring. For ordered
fields~$F$, the completion $\hat F$ is in fact an ordered field, and coincides with the order-theoretic completion
of~$F$ as above.

$\tc$ was originally introduced by Hajnal et al.~\cite{tc0} as a non-uniform complexity class, but following more
recent usage, we define it as the class of languages $L\sset\{0,1\}^*$ recognizable by a $\cxt{DLOGTIME}$-uniform family of
polynomial-size constant-depth circuits using~$\neg$ and unbounded fan-in $\land$, $\lor$, and Majority gates;
equivalently, it consists of languages computable by $O(\log n)$-time threshold Turing machines with $O(1)$ thresholds,
or by constant-time TRAM with polynomially many processors \cite{par-sch}. In terms of descriptive complexity, a
language is in~$\tc$ iff the corresponding class of finite structures is definable in $\cxt{FOM}$, first-order logic
with majority quantifiers~\cite{founif}. A function $F\colon(\{0,1\}^*)^n\to\{0,1\}^*$ is a $\tc$ function if
$\lh{F(X_1,\dots,X_n)}\le p\bigl(\lh{X_1},\dots,\lh{X_n}\bigr)$ for some polynomial $p$, and the bit-graph
$\bigl\{\p{\vec X,i}:\bit\bigl(F(\vec X),i\bigr)=1\bigr\}$ is a $\tc$ predicate. We also consider $\tc$ predicates and
functions where the output or some of the inputs are natural numbers given in unary rather than binary strings; see
\cite[\S IV.3]{cook-ngu} for details. 

We now briefly summarize the definition of~$\vtc$ and its relevant properties, but we refer the reader to
\cite{cook-ngu} (as well as \cite[\S2]{ej:vtcanal}) for more details. $\vtc$ is a theory in a two-sorted first-order
language with equality. The first sort is for natural numbers (called \emph{small} or \emph{unary} numbers), and the
second sort for finite sets of small numbers, which can also be interpreted as binary strings, or as \emph{large} or
\emph{binary} numbers. The second sort is the one we are interested in; the first sort should be thought of as
auxiliary, used for indexing bits of binary numbers. The language of $\vtc$ includes the elementhood predicate~$\in$,
the usual arithmetical functions and predicates $+$, $\cdot$, $0$, $1$, and~$<$ on the first sort, and the $\lh\ $
function whose intended meaning is $\lh X=\sup\{x+1:x\in X\}$. The axioms of~$\vtc$ include several basic axioms
governing the symbols of the language, the comprehension axiom
\[\tag{$\fii\text-\thry{COMP}$} \exists X\le x\:\forall u<x\:\bigl(u\in X\eq\fii(u)\bigr)\]
for $\Sig0$~formulas $\fii$, and an axiom asserting that for any set~$X$, there is a set coding the counting function
$F(i)=\card(X\cap\{0,\dots,i-1\})$ for $i\le\lh X$. Here, the second-order bounded quantifier $\exists X\le x\,\dots$
is defined as $\exists X\,(\lh X\le x\land\dots)$, and similarly for $\forall X\le x\,\dots$; a $\Sig0$~formula has
bounded first-order quantifiers and no second-order quantifiers, and more generally, a $\Sig i$~formula consists of $i$
alternating blocks of second-order bounded quantifiers followed by a $\Sig0$~formula, with the first block being
existential.

All $\tc$ functions have provably total $\Sig1$ definitions in $\vtc$, and $\vtc$ proves comprehension (and therefore
induction over small numbers) for $\Sig0$~formulas in a language expanded with these definable functions, which we will
call \emph{$\tc$~formulas} for short.

$\vtc$ can define (as $\tc$ functions) $+$, $-$, $\cdot$, and $<$ on binary numbers, and proves their basic properties.
(It can also do division with remainder by~\cite{ej:vtcimul}.) If $\sM\model\vtc$, we denote by
$\p{\MN^\sM,+,\cdot,0,1,{<}}$ the second sort of~$\sM$ interpreted as a set of binary natural numbers along with its
arithmetic structure, and extend it with negative numbers to form $\p{\MZ^\sM,+,\cdot,0,1,{<}}$ (the \emph{integers
of~$\sM$}); this is a discretely ordered ring, and in fact, a model of~$\io$ by~\cite{ej:vtc0iopen,ej:vtcimul}. We
define $\p{\MQ^\sM,+,\cdot,0,1,{<}}$ (the \emph{rationals of~$\sM$}) as the fraction field of $\MZ^\sM$, and
$\p{\MR^\sM,+,\cdot,0,1,{<}}$ (the \emph{reals of~$\sM$}) as the completion of~$\MQ^\sM$, which is a real-closed field
by~\cite{ej:vtc0iopen,ej:vtcimul}.

The unary number sort of~$\sM$ embeds (via a $\tc$ function) into $\MN^\sM$ as an initial segment of \emph{logarithmic
numbers}, which we denote $\ML^\sM$. We define the \emph{logarithmically bounded} reals, rationals, and integers by
$\MRL^\sM=\{z\in\MR^\sM:\exists n\in\ML^\sM\:\abs z\le n\}$, $\MQL^\sM=\MQ^\sM\cap\MRL^\sM$, and
$\MZL^\sM=\MZ^\sM\cap\MRL^\sM$. If $n$ is a unary natural number, $2^n$ is represented as a binary number by the set
$\{n\}$. Thus, we can define a $\tc$ function $2^n\colon\ML^\sM\to\MN^\sM$ satisfying $2^1=2$ and $2^{n+m}=2^n2^m$.
(Much more generally, $\vtc$ has a well-behaved definition of products $\prod_{i<n}X_i$ of coded sequences of binary
numbers by~\cite{ej:vtcimul}.)

Let $\texp$ denote the axiom of totality of integer exponentiation; in the context of $\vtc$, it can be simply
expressed as $\ML=\MN$. $\vtc+\texp$ is essentially identical to the common theory $\idze$: if
$\sM\model\vtc+\texp$, the embedding of the unary sort in the binary sort becomes an isomorphism w.r.t.\
$\p{+,\cdot,0,1,{<}}$, and $\MN^\sM\model\idze$. Conversely, a model of $\idze$ expands to a model of $\vtc+\texp$ with
two identical sorts and elementhood predicate defined by $x\in X$ iff $\fl{X/2^x}$ is odd. 

Earlier, Johannsen and Pollett~\cite{joh-pol:d1cr} defined a theory $\dicr$ in the usual one-sorted language of
arithmetic (expanded with a few functions symbols following Buss~\cite{buss}, but these can be in principle eliminated
as they are definable in the $\p{+,\cdot,0,1,{<}}$ language). This theory is bi-interpretable with $\vtc$ such that the
second sort of $\vtc$ becomes the universe of $\dicr$, hence models of $\dicr$ are exactly the structures $\MN^\sM$ for
$\sM\model\vtc$. Thus, we could have formulated everything more directly in terms of models of $\dicr$;
nevertheless, we use $\vtc$ as it became a de facto standard theory corresponding to~$\tc$.

Recursive saturation was introduced by Barwise and Schlipf~\cite{bar-sch}. Let $\sM=\p{M,\dots}$ be a
structure in a finite language~$\lang$. If $\vec a\in M$ and $\Gamma(x,\vec y)$ is a recursive set of $\lang$-formulas,
then $\Gamma(x,\vec a)$ is a \emph{recursive type} of~$\sM$, which is \emph{finitely satisfiable} if
$\sM\model\exists x\,\ET_{\fii\in\Gamma'}\fii(x,\vec a)$ for each finite $\Gamma'\sset\Gamma$, and
\emph{realized} by $c\in M$ if $\sM\model\Gamma(c,\vec a)$. Then $\sM$ is \emph{recursively saturated} if
every finitely satisfiable recursive type of~$\sM$ is realized in~$\sM$. By Craig's trick, this definition does not
change if we consider recursively enumerable types or $\tc$~types in place of recursive types.

Two structures $\mathfrak A$ and~$\mathfrak B$ are \emph{jointly recursively saturated} if a structure
$\p{\mathfrak A,\mathfrak B}$ encompassing both in a suitable way is recursively saturated. The uniqueness theorem
states that elementarily equivalent countable jointly recursively saturated structures are isomorphic. We will not work
with $\p{\mathfrak A,\mathfrak B}$ as such, but in view of the fact that recursive saturation is preserved by
interpretation, we can rephrase the uniqueness theorem as follows:
\begin{Thm}[Barwise and Schlipf~\cite{bar-sch}]\label{thm:rec-sat}
Let $\mathfrak A$ and~$\mathfrak B$ be elementarily equivalent countable structures interpretable in a recursively
saturated structure~$\sM$. Then $\mathfrak A\simeq\mathfrak B$.
\noproof\end{Thm}

\section{Real exponential in models of $\vtc$}\label{sec:real-expon-models}

Let $\sM$ be a model of~$\vtc$. How can we show that $\sM$ (more precisely, the discretely ordered ring $\MZ^\sM$) is
an exponential integer part of a real-closed exponential field? Our starting points are the result of
\cite{ej:vtc0iopen,ej:vtcimul} that $\MR^\sM$ is a real-closed field with integer part~$\MZ^\sM$, and the construction
of a natural exponential function on~$\MR^\sM$ in~\cite{ej:vtcanal}. Stated for base-$2$ exponentiation, the relevant
properties of the latter can be summarized as follows:
\begin{Thm}[Je\v r\'abek \cite{ej:vtcanal}]\label{thm:vtcanal}
For any model $\sM\model\vtc$, the usual function $2^n\colon\ML^\sM\to\MN^\sM$ extends to an ordered group isomorphism
$2^x\colon\p{\MRL^\sM,+,0,{<}}\to\p{\MR_{>0}^\sM,\cdot,1,{<}}$.
\noproof\end{Thm}

This almost shows that $\MR^\sM$ is a real-closed exponential field, and $\MZ^\sM$ is its exponential IP, were it not
for the pesky ${}_\ML$ in the domain of~$2^x$. Which is, of course, essential: if integer exponentiation is not total,
we can only expect a reasonably well-behaved real exponential function to be defined on logarithmically small numbers,
and even if we manage to find a wild exponential defined on all of~$\MR^\sM$, there is no way it could be compatible
with the usual integer $2^n$ function. Let us state for the record that we are done if integer exponentiation \emph{is}
total, though (recall that $\vtc+\texp=\idze$):
\begin{Cor}\label{cor:idz+exp}
Any model $\sM\model\idze$ is an exponential IP of a real-closed exponential field
$\p{\MR^\sM,+,\cdot,0,1,{<},2^x}$.
\noproof\end{Cor}

But Theorem~\ref{thm:vtcanal} makes significant progress even if $\sM\nmodel\texp$:
\begin{Cor}\label{cor:r-rl-eip}
Let $\sM\model\vtc$, and assume there exists an isomorphism
\[f\colon\p{\MR^\sM,\MZ^\sM,+,0,1,{<}}\to\p{\MRL^\sM,\MZL^\sM,+,0,1,{<}}\]
such that $2^{f(x)}>x$ for all $x\in\MR^\sM_{>0}$. Then $\exp(x)=2^{f(x)}$ makes $\MR^\sM$ a real-closed exponential
field with exponential IP $\MZ^\sM$.
\noproof\end{Cor}

Our basic idea is to construct such an isomorphism~$f$ using Theorem~\ref{thm:rec-sat}. With any luck,
$\p{\MR^\sM,\MZ^\sM,+,0,1,{<}}$ and $\p{\MRL^\sM,\MZL^\sM,+,0,1,{<}}$ will be elementarily equivalent. However, these
structures are uncountable even if $\sM$ itself is countable, hence we cannot directly apply Theorem~\ref{thm:rec-sat} to
them; moreover, they are not interpretable in~$\sM$, which leads to difficulties when trying to establish they are
jointly recursively saturated.

One way to get around these problems is to use the fact that any isomorphism of ordered groups extends to an
isomorphism of their completions, thus it is enough to construct an isomorphism
$f\colon\p{\MQ^\sM,\MZ^\sM,+,0,1,{<}}\to\p{\MQL^\sM,\MZL^\sM,+,0,1,{<}}$. Then the original strategy essentially works:
using a quantifier elimination result, we can prove that $\p{\MQ^\sM,\MZ^\sM,+,0,1,{<}}$ and
$\p{\MQL^\sM,\MZL^\sM,+,0,1,{<}}$ are elementarily equivalent, and $\p{\MQ^\sM,\MZ^\sM,\MQL^\sM,+,0,1,{<}}$ is
recursively saturated. If $\sM$ is countable, this implies that $\p{\MQ^\sM,\MZ^\sM,+,0,1,{<}}$ and
$\p{\MQL^\sM,\MZL^\sM,+,0,1,{<}}$ are isomorphic. This argument per se does not ensure the growth condition
$2^{f(x)}>x$, which requires yet more work.

An even easier route is to use the following observation:
\begin{Lem}\label{lem:z-r}
If $\sM\model\vtc$, any isomorphism $f\colon\p{\MZ^\sM,+,0,1,{<}}\to\p{\MZL^\sM,+,0,1,{<}}$ extends to an isomorphism
$\ob f\colon\p{\MR^\sM,\MZ^\sM,+,0,1,{<}}\to\p{\MRL^\sM,\MZL^\sM,+,0,1,{<}}$.
\end{Lem}
\begin{Pf}
Define $\ob f(x)=f\bigl(\fl x\bigr)+\{x\}$. It is clear that $\ob f$ is an order-preserving bijection, thus we only
need to check that it is a group homomorphism. Given $x,y\in\MR^\sM$, either $\{x\}+\{y\}\in[0,1)$ or
$\{x\}+\{y\}\in[1,2)$. In the latter case, $\{x+y\}=\{x\}+\{y\}-1$ and $\fl{x+y}=\fl x+\fl y+1$, thus
\begin{align*}
\ob f(x+y)&=f\bigl(\fl x+\fl y+1\bigr)+\{x\}+\{y\}-1\\
&=f\bigl(\fl x\bigr)+f\bigl(\fl y\bigr)+f(1)+\{x\}+\{y\}-1=\ob f(x)+\ob f(y)
\end{align*}
as $f(1)=1$. In the former case, $\fl{x+y}=\fl x+\fl y$ and $\{x+y\}=\{x\}+\{y\}$, thus $\ob f(x+y)=\ob f(x)+\ob f(y)$
by a similar (easier) argument.
\end{Pf}

Consequently, we could make do with recursive saturation of just $\p{\MZ^\sM,\MZL^\sM,+,0,1,{<}}$. However, we consider
the recursive saturation result to be of independent interest in its own right, and therefore proceed to prove it in
full generality for $\p{\MQ^\sM,\MZ^\sM,\MQL^\sM,+,0,1,{<}}$ as suggested above.

\section{The theory of three groups}\label{sec:theory-three-groups}

Our first task is to axiomatize the theory of $\p{\MQ^\sM,\MZ^\sM,\MQL^\sM,+,0,1,{<}}$ and show that it enjoys
quantifier elimination down to a convenient class of formulas.

\begin{Def}\label{def:3g}
The \emph{theory of three groups} (denoted $\GGG$) is a first-order theory in the language $\lang_\GGG=\p{Z,L,+,0,1,{<}}$,
where $Z$ and~$L$ are unary predicates. We will often treat $Z$ and~$L$ as sets, writing $x\in Z$ for $Z(x)$, and using
it as quantifier bounds such as $\exists x\in Z\,\dots$; we will also denote the whole universe as~$Q$. The axioms of
$\GGG$ are:
\begin{enumerate}
\item\label{item:1} $\p{Q,+,0,{<}}$ is a divisible ordered group.
\item\label{item:2} $Z$ is an integer part of~$Q$ with a least positive element~$1$.
\item\label{item:3} $L$ is a convex subgroup of~$Q$ containing~$1$.
\end{enumerate}
Notice that the axioms imply that $Z$ is a $\Z$-group. We define $qx$ for $q\in\Q$ and $x\in Q$ as usual (being a
torsion-free divisible group, $Q$ carries a definable structure of a $\Q$-linear space), and write $q1$ as just~$q$; we
also write $x\equiv y\pmod m$ for $x-y\in mZ$, where $m\in\N_{>0}$.
\end{Def}

\begin{Exm}\label{exm:vtc-3g}
For any $\sM\model\vtc$, $\p{\MQ^\sM,\MZ^\sM,\MQL^\sM,+,0,1,{<}}$ and $\p{\MR^\sM,\MZ^\sM,\MRL^\sM,+,0,1,{<}}$ are models
of~$\GGG$.
\end{Exm}

Notice that $\GGG$ is incomplete, as it does not decide the sentence $Q=L$ (i.e., $\forall x\,L(x)$).

\begin{Def}\label{def:special-fla}
A \emph{special formula} is a Boolean combination of formulas of the form
\begin{align}
\label{eq:1}\tsum_in_i\{x_i\}&\ge n,\\
\label{eq:2}\tsum_in_i\fl{x_i}&\ge n,\\
\label{eq:3}\fl{x_i}&\equiv k\pmod m,\\
\label{eq:4}\tsum_in_ix_i&\in L,\\
\label{eq:5}Q&=L,
\end{align}
where $n_i,n,k,m\in\Z$, $0\le k<m$.
\end{Def}

\begin{Thm}\label{thm:qe-3g}
In $\GGG$, every formula is equivalent to a special formula.
\end{Thm}
\begin{Pf}
First, any formula is equivalent to one where the only atomic formulas are of the form $x+y=z$, $x\ge0$,
$x=1$, $Z(x)$, or $L(x)$ for some variables $x,y,z$. These are easy to express by special formulas: e.g., $x+y=z$ is
equivalent to
\[\bigl(\{x\}+\{y\}=\{z\}\land\fl x+\fl y=\fl z\bigr)\lor\bigl(\{x\}+\{y\}=\{z\}+1\land\fl x+\fl y=\fl z-1\bigr),\]
which can be further rewritten in terms of inequalities. Thus, it suffices to show that special formulas are closed
under existential quantification up to equivalence.

Let us consider a formula $\fii(\vec x)=\exists x\,\theta(x,x_0,\dots,x_{t-1})$, where $\theta$ is special. Using
standard manipulations (replacing negated inequalities and congruences, writing $\theta$ in DNF, commuting $\lor$
with~$\exists$, moving out conjuncts without~$x$), we may assume $\theta(x,\vec x)=\ET_j\theta_j(x,\vec x)$, where each
$\theta_j$ has the form
\begin{align}
\label{eq:6}n\{x\}&=\ell(\{\vec x\}),\\
\label{eq:7}n\{x\}&>\ell(\{\vec x\}),\\
\label{eq:8}n\fl x&\ge\ell(\fl{\vec x}),\\
\label{eq:9}\fl x&\equiv k\pmod m,\\
\label{eq:10}nx-\ell(\vec x)&\in L.
\end{align}
Here, $n\in\Z\bez\{0\}$, $0\le k<m$, $\ell(\vec x)=\sum_{i<t}n_ix_i+r$ with $n_i,r\in\Z$, $\{\vec x\}$ denotes
$\{x_0\},\dots,\{x_{t-1}\}$, and similarly for $\fl{\vec x}$. (Since $u>v\iff u\ge v+1$ for $u,v\in Z$, we do not need
a version of~\eqref{eq:8} with strict inequality.) Notice that \eqref{eq:10} is equivalent to
\[n\fl x-\ell(\fl{\vec x})\in L,\]
thus we can write $\theta$ in the form $\theta'(\{x\},\{\vec x\})\land\theta''(\fl x,\fl{\vec x})$. For every $u\in Z$
and $v\in[0,1)$, there is $x$ such that $\{x\}=u$ and $\fl x=v$, namely $x=u+v$; it follows that $\exists
x\,\theta(x,\vec x)$ is equivalent to
\[\exists x\in[0,1)\:\theta'(x,\{\vec x\})\land\exists x\in Z\:\theta''(x,\fl{\vec x}),\]
where $\theta'$ is a conjunction of formulas of the form \eqref{eq:6} and~\eqref{eq:7}, and $\theta''$ is a conjunction
of formulas of the form \eqref{eq:8}, \eqref{eq:9}, and~\eqref{eq:10}. The first part can be further rewritten as 
\[\theta'(0,\{\vec x\})\lor\exists x\:\bigl(x>0\land x<1\land\theta'(x,\{\vec x\})\bigr);\]
the first disjunct is a special formula, hence we may ignore it, and then we may just assume that $\theta'$ includes
$x>0$ and $x<1$ among the inequalities~\eqref{eq:7}. If $\theta'$ includes any equality~\eqref{eq:6}, then
$\exists x\:\theta'(x,\{\vec x\})$ is equivalent to $\theta'\bigl(\frac1n\ell(\{\vec x\}),\{\vec x\}\bigr)$, which is
a conjunction of linear%
\footnote{We allow linear functions, equations, and inequalities to be inhomogeneous, i.e., of the form $\sum_iq_ix_i+q$.}
equations and inequalities in $\{\vec x\}$ with rational coefficients; multiplying each
(in)equality by $\abs n$, we obtain a special formula. Otherwise, $\theta'$ consists only of strict
inequalities~\eqref{eq:7}. Dividing each inequality by the coefficient of~$x$, we can write $\theta'$ as
\[\ET_{i\in I^+}x<\ell_i(\{\vec x\})\land\ET_{i\in I^-}x>\ell_i(\{\vec x\}),\]
where $\ell_i$ are linear functions with rational coefficients. Then $\exists x\,\theta'(x,\{\vec x\})$ is
equivalent to
\[\ET_{i^+\in I^+}\ET_{i^-\in I^-}\ell_{i^+}(\{\vec x\})>\ell_{i^-}(\{\vec x\}),\]
which can be written as a special formula.

It remains to deal with $\exists x\in Z\:\theta''(x,\fl{\vec x})$. In order to simplify the notation, we will assume
$\vec x$ are given as elements of~$Z$ so that we can henceforth drop the $\fl{\dots}$ signs. Multiplying the
inequalities \eqref{eq:8} and the expressions in~\eqref{eq:10} by suitable constants, we can ensure that they all use
the same $n$ up to sign. Replacing also $x\equiv k\pmod m$ with $nx\equiv nk\pmod{nm}$, we can then write
$\theta''(x,\vec x)$ so that $x$ occurs everywhere with a multiplier $\pm n$. Using
\[\exists x\in Z\:\psi(nx,\vec x)\iff\exists x\in Z\:\bigl(\psi(x,\vec x)\land x\equiv0\pmod n\bigr),\]
we reduce the problem to the case $n=1$.

Moreover, we can combine the congruences \eqref{eq:9} using the Chinese remainder theorem: the conjunction of
$x\equiv k_0\pmod{m_0}$ and $x\equiv k_1\pmod{m_1}$ is equivalent either to~$\bot$, if
$k_0\nequiv k_1\pmod{\gcd(m_0,m_1)}$, or to $x\equiv k\pmod m$, where $m=\lcm(m_0,m_1)$ and $k\equiv k_i\pmod{m_i}$.
Thus, we can write $\theta''$ as
\[x\equiv k\pmod m\land\ET_{i\in I^+}x\le\ell_i(\vec x)\land\ET_{i\in I^-}x\ge\ell_i(\vec x)\land
  \ET_{i\in J^+}x-\ell_i(\vec x)\in L\land\ET_{i\in J^-}x-\ell_i(\vec x)\notin L,\]
where $\ell_i$ are linear functions with integer coefficients. We may also assume $I^+\cup I^-=J^+\cup J^-$: if, say,
$i\in J^+\cup J^-\bez(I^+\cup I^-)$, we have
\[\exists x\in Z\:\theta''(x,\vec x)\iff\exists x\in Z\:\bigl(\theta''(x,\vec x)\land x\ge\ell_i(\vec x)\bigr)
  \lor\exists x\in Z\:\bigl(\theta''(x,\vec x)\land x\le\ell_i(\vec x)\bigr),\]
and likewise for $i\in I^+\cup I^-\bez(J^+\cup J^-)$.

Assume first $J^+=\nul$, thus $J^-=I^+\cup I^-$. We claim that $\exists x\in Z\,\theta''(x,\vec x)$ is equivalent to
\begin{equation}\label{eq:11}
\ET_{i^+\in I^+}\ET_{i^-\in I^-}\bigl(\ell_{i^+}(\vec x)\ge\ell_{i^-}(\vec x)
  \land\ell_{i^+}(\vec x)-\ell_{i^-}(\vec x)\notin L\bigr),
\end{equation}
which can be written as a special formula. It is easy to see that $\exists x\in Z\,\theta''(x,\vec x)$
implies~\eqref{eq:11}. For the converse, \eqref{eq:11} expresses that if $\ell^-=\max\{\ell_i(\vec x):i\in I^-\}$
and $\ell^+=\min\{\ell_i(\vec x):i\in I^+\}$, then $\ell^+>\ell^-$ and $\ell^+-\ell^-\notin L$. Splitting the interval
$[\ell^-,\ell^+]$ in thirds, we can find $\ell^-<u<v<\ell^+$ such that $u-\ell^-$, $v-u$, and $\ell^+-v$ are still
outside~$L$; in particular, $v-u$ is infinite, hence there is $u<x<v$ such that $x\equiv k\pmod m$. Then for each
$i\in I^+\cup I^-$, $x-\ell_i(\vec x)\notin L$ and has the right sign, thus $\theta''(x,\vec x)$.

This discussion tacitly assumed $I^+,I^-\ne\nul$. If $I^+=I^-=\nul$, $\exists x\in Z\,\theta''(x,\vec x)$ is
always true, as is~\eqref{eq:11} (vacuously). However, if $I^+=\nul\ne I^-$, we need to assume $Q\ne L$ to find
$u>\ell^-$ such that $u-\ell^-\notin L$; on the other hand, $\exists x\in Z\,\theta''(x,\vec x)$ clearly implies
$Q\ne L$ as $J^-\ne\nul$. Thus, if $I^+=\nul\ne I^-$ or $I^+\ne\nul=I^-$, then $\exists x\in Z\,\theta''(x,\vec x)$
is equivalent to the special formula $Q\ne L$ rather than to~\eqref{eq:11}.

Finally, assume $J^+\ne\nul$. Substituting $x+\ell_i(\vec x)$ for $x$ if necessary (which does not change the truth
value of $\exists x\in Z\,\theta''(x,\vec x)$), we may assume $\theta''$ includes a conjunct $x\in L$. But on
condition of $x\in L$, we can dispense with the remaining conjuncts involving $L$, as
$x-\ell_i(\vec x)\in L$ iff $\ell_i(\vec x)\in L$, which can be moved outside the scope of the $\exists x$ quantifier. Thus, $\theta''$ simplifies to
\[x\in L\land x\equiv k\pmod m\land\ET_{i\in I^+}x\le\ell_i(\vec x)\land\ET_{i\in I^-}x\ge\ell_i(\vec x).\]
We claim that $\exists x\in Z\,\theta''(x,\vec x)$ is equivalent to
\begin{align*}
\ET_{i\in I^+}\bigl(\ell_i(\vec x)\ge0\lor\ell_i(\vec x)\in L\bigr)
&\land\ET_{i\in I^-}\bigl(\ell_i(\vec x)\le0\lor\ell_i(\vec x)\in L\bigr)\\
&\land\ET_{i_+\in I^+}\ET_{i^-\in I^-}\exists x\in Z\:
    \bigl(\ell_{i^-}(\vec x)\le x\le\ell_{i^+}(\vec x)\land x\equiv k\pmod m\bigr).
\end{align*}
If this formula holds, let $\ell_-$ and $\ell_+$ be as above (assuming $I^+,I^-\ne\nul$). The first two conjuncts
ensure that $[\ell_-,\ell_+]$ intersects~$L$. If the convex set $[\ell_-,\ell_+]\cap L$ has length at least~$m$, it
contains an $x\equiv k\pmod m$, which witnesses $\theta''(x,\vec x)$. Otherwise we must have $\ell_-,\ell_+\in L$, and
the third conjunct ensures there is $x\in[\ell_-,\ell_+]$ such that $x\equiv k\pmod m$, which then belongs to~$L$ as
well. It is easy to see that the equivalence holds even if $I^+$ or $I^-$ is empty.

Each of the formulas $\exists x\in Z\,\bigl(\ell_{i^-}(\vec x)\le x\le\ell_{i^+}(\vec x)\land x\equiv k\pmod m\bigr)$
is equivalent to
\[\LOR_{\substack{0\le\vec a,a<m\\\ell_{i^+}(\vec a)\equiv a+k\pmod m}}
  \Bigl(\ET_{j<t}x_j\equiv a_j\pmod m\land\ell_{i^-}(\vec x)\le\ell_{i^+}(\vec x)-a\Bigr),\]
as the reader can check.
\end{Pf}
\begin{Cor}\label{cor:3g-compl}
The only completions of $\GGG$ are $\GGG+Q=L$ and $\GGG+Q\ne L$.
\noproof\end{Cor}
\begin{Def}\label{def:2g}
Let $\GG$ denote the theory in the language $\lang_\GG=\p{Z,+,0,1,{<}}$ axiomatized by \ref{item:1} and~\ref{item:2} from
Definition~\ref{def:3g}.
\end{Def}
\begin{Cor}\label{cor:2g}
The theory $\GG$ is complete. Any formula is in $\GG$ equivalent to a Boolean combination of formulas of the form
\eqref{eq:1}, \eqref{eq:2}, and~\eqref{eq:3}.
\end{Cor}
\begin{Pf}
$\GG$ is essentially identical to $\GGG+Q=L$.
\end{Pf}

\section{Recursive saturation of models of $\GGG$}\label{sec:recurs-satur-models}

Our goal is to show that $\GGG$ reducts of nonstandard models $\sM\model\vtc$ are recursively saturated. The key
ingredient of the proof will be a $\tc$ truth predicate for (a subset of) special formulas; this idea works nicely for
$\lang_\GG$-formulas of the form \eqref{eq:1}--\eqref{eq:3}, but fails miserably for formulas of the form \eqref{eq:4},
as $\ML^\sM$ is not definable in~$\sM$ by any bounded formula (unless $\sM\model\texp$). To get around this
problem, we give in this section a description of recursive saturation of models of $\GGG$ that separates the roles of
$\lang_\GG$ and~$L$.

\begin{Def}\label{def:cof}
If $\p{P,\le}$ is a poset and $X\sset P$, we define $X^\up=\{u\in P:X\le u\}$, where as usual, $X\le u$ means $\forall
x\,(x\in X\to x\le u)$. A subset $Y\sset X$ is \emph{(upwards) cofinal in $X$} if $\forall x\in X\,\exists y\in Y\,x\le
y$, and it is \emph{downwards cofinal in $X$} if $\forall x\in X\,\exists y\in Y\,y\le x$. If $Y\sset P$ is downwards
cofinal in $X^\up$, we also say that $Y$ is \emph{cofinal above $X$}.

Let $G$ be a divisible ordered group, hence a $\Q$-linear space. We write $Sa=\{qa:q\in S\}$ for any
$a\in G$ and $S\sset\Q$ such as $S=\N$ or $S=\N^{-1}$, the latter abbreviating $\{n^{-1}:n\in\N_{>0}\}$. If $X\sset G$,
let $\p X_\Q$ denote the $\Q$-linear span of~$X$.
\end{Def}

Observe that if $\p{Q,Z,L,+,0,1,{<}}$ is a recursively saturated model of $\GGG+Q\ne L$, then no set of the form $\N a$
is cofinal in $L$ as $\Gamma(x)=\{x>na:n\in\N\}$ is a recursive type, and likewise, no set of the form $\N^{-1}a$ is
cofinal above~$L$. We could generalize this observation to finitely generated subsets in place of $\N a$ or $\N^{-1}a$,
but this would be equivalent due to the following lemma.
\begin{Lem}\label{lem:fg-cof}
Let $G$ be a divisible ordered group, finite-dimensional as a $\Q$-linear space, and $C\sset G$ a proper convex
subgroup.
\begin{enumerate}
\item\label{item:4} There exists $a\in G$ such that $\N a$ is cofinal in~$C$.
\item\label{item:5} There exists $a\in G$ such that $\N^{-1}a$ is cofinal above~$C$.
\end{enumerate}
\end{Lem}
\begin{Pf}

\ref{item:4}: $C$ is a linear subspace of~$G$, thus $C=\p{a_i:i<r}_\Q$ for some $\vec a\in C$. Putting
$a=\max_i\abs{a_i}$, $\N a$ is cofinal in~$C$, as $\sum_iq_ia_i\le a\sum_i\CL[big]{\abs{q_i}}$ for all $\vec q\in\Q$.

\ref{item:5}: There are only finitely many convex subgroups of~$G$, as they form a family of linear subspaces totally
ordered by inclusion. Thus, there exists a minimal convex subgroup $C'\sSset C$; then $\N^{-1}a$ is cofinal above~$C$
for any $a\in C'_{>0}\bez C$, as $C''=\{x\in G:\forall n\in\N_{>0}\,\abs x\le n^{-1}a\}$ is a convex subgroup of~$G$
such that $C\sset C''\ssset C'$.
\end{Pf}

We now prove a characterization of recursive saturation of models of $\GGG$. In view of the discussion above, it shows
that obvious necessary conditions are also sufficient.
\begin{Thm}\label{thm:3g-rec-sat}
A model $\sM=\p{Q,Z,L,+,0,1,{<}}\model\GGG$ is recursively saturated if and only if
\begin{enumerate}
\item\label{item:6} $\sM\res\lang_\GG$ is recursively saturated, and
\item\label{item:7} there is no $a\in Q$ such that $\N a$ is cofinal in~$L$ or $\N^{-1}a$ is cofinal
above~$L$.
\end{enumerate}
\end{Thm}
\begin{Pf}
We have already seen that if $\sM$ is recursively saturated, it satisfies \ref{item:6} and~\ref{item:7}. Conversely,
assume that \ref{item:6} and~\ref{item:7} hold, and let $\Gamma(x,\vec a)$ be a finitely satisfiable recursive type; we
have to show that $\Gamma(x,\vec a)$ is realized in~$\sM$. We may assume $Q\ne L$, and $\vec a\sset Z\cup(0,1)$, which
ensures that each $\fl{a_i}$ or $\{a_i\}$ is either $a_i$ itself or~$0$. We also assume $\vec a$ includes~$1$. By Theorem~\ref{thm:qe-3g}, we may assume that
$\Gamma$ consists of special formulas (not involving~\eqref{eq:5}, as this can be replaced with~$\bot$); by eliminating
formulas~\eqref{eq:4} in a suitable way, we will construct a finitely satisfiable recursive $\lang_\GG$-type that
implies~$\Gamma$, and appeal to~\ref{item:6}. We distinguish two cases.

\paragraph{Case 1:} $\Gamma(x,\vec a)\cup\{x-\alpha\in L\}$ is finitely satisfiable for some $\alpha\in\p{\vec a}_\Q$.
Substituting $x+\alpha$ for~$x$ if necessary, we may assume $\alpha=0$. By Lemma~\ref{lem:fg-cof}, $\p{\vec a}_\Q\cap L$
has a cofinal subset of the form $\N a$, hence it is not cofinal in~$L$ due to~\ref{item:7}; thus, we may fix $b\in L$
such that $\p{\vec a}_\Q\cap L<b$. Let $\Gamma'(x,\vec a,b)$ be the $\lang_\GG$-type obtained from $\Gamma$ by
replacing each subformula of the form $nx+\ell(\vec a)\in L$ (where $\ell$ is a $\Z$-linear function) with
$\abs{\ell(\vec a)}<b$. Then $\Gamma(x,\vec a)\cup\{x\in L\}$ is equivalent to $\Gamma'(x,\vec a,b)\cup\{x\in L\}$, as
$x\in L$ implies
\[nx+\ell(\vec a)\in L\iff\ell(\vec a)\in L\iff\abs{\ell(\vec a)}<b.\]
Put $\Gamma''(x,\vec a,b)=\Gamma'(x,\vec a,b)\cup\{\abs x<b\}$. Since $(-b,b)\sset L$, any realizer of the
$\lang_\GG$-type $\Gamma''$ also realizes~$\Gamma$. It remains to show that $\Gamma''$ is finitely satisfiable, whence
realized by~\ref{item:6}.

Let $\fii$ be the conjunction of a finite subset of $\Gamma'(x,\vec a,b)$; we need to satisfy $\fii$ by an element of
$(-b,b)$. Since $\Gamma'\cup\{x\in L\}$ is finitely satisfiable, there exists $u\in L$ such that
$\sM\model\fii(u,\vec a,b)$. Notice that $b$ only occurs in~$\fii$ in subformulas of the form $\abs{\ell(\vec a)}<b$,
each of which has a fixed truth value independent of~$x$, and as such can be eliminated; thus, we may assume that $b$
does not occur in~$\fii$, i.e., $\fii$ is a special $\lang_\GG$-formula in $x$ and~$\vec a$. By writing $\fii$ in DNF
and separating $x$ to one side, $u$ satisfies in~$\sM$ a conjunction of formulas of the form
\[\{x\}\leqg\ell(\vec a),\qquad
\fl{x}\leqg\ell(\vec a),\qquad
\fl{x}\equiv k\pmod m\]
that implies $\fii(x,\vec a)$, where $\ell$ are $\Q$-linear functions and ${\leqg}\in\{{<},{=},{>}\}$. That
is, there exist (possibly degenerate) intervals $I\sset(0,1)$ and $J$ with endpoints in $\p{\vec
a}_\Q\cup\{\pm\infty\}$, and an arithmetic progression $P\sset Z$ with standard modulus~$m$, such that
\begin{equation}\label{eq:12}
u\in I+(J\cap P)\sset\{x\in Q:\sM\model\fii(x,\vec a)\}.
\end{equation}
Shortening $J$ and/or negating $x$ if necessary, we may assume $J=(v,w)$ with $0\le v<u<w\le+\infty$. Since $u\in L$, we
have $v\in\p{\vec a}_\Q\cap L$, hence $v<b$. If $w\in L$ as well, we have $u<w<b$ and we are done; otherwise $w-v$
is nonstandard, and $v+t\in J\cap P$ for some $0<t\le m$, thus $\fii$ is satisfied by $v+t+\{u\}<b$.

\paragraph{Case 2:} For all $\alpha\in\p{\vec a}_\Q$, $\Gamma(x,\vec a)\cup\{x-\alpha\in L\}$ is finitely
unsatisfiable, i.e., $\Gamma(x,\vec a)$ implies $x-\alpha\notin L$; thus, $\Gamma(x,\vec a)\cup(x\notin\p{\vec a}_\Q+L)$
is finitely satisfiable, where $(x\notin\p{\vec a}_\Q+L)$ is the type
\[\bigl\{x-\ell(\vec a)\notin L:\ell\text{ is a $\Q$-linear function}\bigr\}.\]
Using Lemma~\ref{lem:fg-cof}, we can fix $b\in L^\up$ such that $b<\p{\vec a}_\Q\cap L^\up$. Let $\Gamma'(x,\vec a,b)$ be
obtained from $\Gamma(x,\vec a)$ by replacing each subformula of the form $nx+\ell(\vec a)\in L$ with $\bot$ if
$n\ne0$, and with $\abs{\ell(\vec a)}<b$ if $n=0$. Notice that $\Gamma(x,\vec a)\cup(x\notin\p{\vec a}_\Q+L)$ and
$\Gamma'(x,\vec a,b)\cup(x\notin\p{\vec a}_\Q+L)$ are equivalent. Put
\[\Gamma''(x,\vec a,b)=
   \Gamma'(x,\vec a,b)\cup\bigl\{\abs{x-\ell(\vec a)}>b:\ell\text{ is a $\Q$-linear function}\bigr\}.\]
Again, the $\lang_\GG$-type $\Gamma''(x,\vec a,b)$ implies $\Gamma(x,\vec a)$, hence we only need to show that
$\Gamma''(x,\vec a,b)$ is finitely satisfiable to finish the proof.

Let $\fii$ be the conjunction of a finite subset of $\Gamma'(x,\vec a,b)$, and $A$ a finite subset of~$\p{\vec a}_\Q$;
we will satisfy the formula $\fii(x,\vec a,b)\land\ET_{\alpha\in A}\abs{x-\alpha}>b$ in~$\sM$. As in Case~1, we may
assume that $b$ does not occur in~$\fii$ by eliminating subformulas that do not depend on~$x$. We may also assume that
$\frac1n\ell(\vec a)\in A$ for every subformula $n\fl x\leqg\ell(\vec a)$ that occurs in~$\fii$. Since
$\Gamma'(x,\vec a,b)\cup(x\notin\p{\vec a}_\Q+L)$ is finitely satisfiable, there exists $u\in Q$ such that
$\sM\model\fii(u,\vec a,b)$ and $u-\alpha\notin L$ for all $\alpha\in A$. As above, there exist (possibly degenerate)
intervals $I\sset(0,1)$ and $J$ with endpoints in $\p{\vec a}_\Q\cup\{\pm\infty\}$, and an arithmetic progression
$P\sset Z$ with standard modulus~$m$, such that \eqref{eq:12} holds; moreover, the endpoints of~$J$, if finite, belong
to~$A$, thus by possibly shortening~$J$, we may assume $J=(\alpha_-,\alpha_+)$, where
\begin{align*}
\alpha_+&=\rlap{$\min$}\phantom\max\bigl\{\alpha\in A\cup\{+\infty\}:\alpha>u\bigr\},\\
\alpha_-&=\max\bigl\{\alpha\in A\cup\{-\infty\}:\alpha<u\bigr\}.
\end{align*}
Since $\alpha_+-u,u-\alpha_-\notin L$, we see that
$\frac13(\alpha_+-\alpha_-)\in L^\up\cap\p{\vec a}_\Q\cup\{+\infty\}$, hence $\frac13(\alpha_+-\alpha_-)>b$. Splitting
$J$ in thirds, the middle part contains an element $u'$ such that $\fl{u'}\in P$ and $\{u'\}\in I$; then
$\sM\model\fii(u',\vec a)$ and $\alpha_+-u',u'-\alpha_->b$, thus $\abs{u'-\alpha}>b$ for all $\alpha\in A$.
\end{Pf}

Given that every $\lang_\GG$-formula is a Boolean combination of formulas of $\p{Z,+,0,1,{<}}$ and formulas of
$\p{Q,+,0,1,{<}}$ (restricted to $[0,1]$), one may wonder whether recursive saturation of models of~$\GG$ can be
further characterized in terms of recursive saturation of the $\p{Z,+,0,1,{<}}$ and $\p{Q,+,0,1,{<}}$ reducts. However,
it is not as simple as that. For example, we have a notion of a ``standard system'' for either reduct:
a set $X\sset\N$ can be encoded by binary expansion of an $x_Z\in Z$ (i.e., $n\in X$ iff
$\LOR_{a<2^n}x_Z\equiv 2^n+a\pmod{2^{n+1}}$), or by binary expansion of an $x_Q\in[0,1)$ (i.e., $n\in X$ iff
$\LOR_{a<2^n}2a+1\le2^{n+1}x_Q<2a+2$). If $\p{Q,Z,+,0,1,{<}}$ is recursively saturated, any coinfinite set
represented in $Z$ is represented in $(0,1)$, and vice versa, so the two reducts interact in a nontrivial way.

\section{$\GGG$ reducts of models of $\vtc$}\label{sec:ggg-reducts-models}

We are now ready to prove our main results.

\begin{Thm}\label{thm:vtc-rec-sat}
If $\sM\model\vtc$ is nonstandard, then $\p{\MQ^\sM,\MZ^\sM,\MQL^\sM,+,0,1,{<}}$ is recursively saturated.
\end{Thm}
\begin{Pf}
Let $c\in\ML^\sM\bez\N$. For each $a\in\MQ_{\ML,>0}^\sM$, $\N a<ca\in\MQL^\sM$, hence $\N a$ is not cofinal in $\MQL^\sM$, and
for each $a\in\MQ_{>0}^\sM\bez\MQL^\sM$ (if any), $\MQL^\sM<c^{-1}a<\N^{-1}a$, hence $\N^{-1}a$ is not cofinal
above~$\MQL^\sM$. Thus, in view of Theorem~\ref{thm:3g-rec-sat}, it suffices to prove that $\p{\MQ^\sM,\MZ^\sM,+,{<}}$ is
recursively saturated.

Let $\Gamma(x,\vec a)$ be a finitely satisfiable recursive type, where $\vec a\sset\MZ^\sM\cup(0,1)$ is $\Q$-linearly
independent, and one of the $a_i$ is~$1$. We first syntactically simplify the type while keeping it recursive. By
Corollary~\ref{cor:2g}, we may assume each $\fii\in\Gamma$ to be a special $\lang_\GG$-formula; we may write it in CNF and
split the conjunctions to make each formula a disjunction of formulas of the form \eqref{eq:1}--\eqref{eq:3} or their
negations (with no constant coefficients, and using just $\vec a$ in place of $\{\vec x\}$ and $\fl{\vec x}$, as each
$\{a_i\}$ and $\fl{a_i}$ equals either $a_i$ or~$0$). Negations of \eqref{eq:2} or~\eqref{eq:3} can be replaced with
(disjunctions of) positive formulas of the same type. Formulas \eqref{eq:1} and their negations can be written as
disjunctions of strict inequalities and equations. The latter can be eliminated: if $\Gamma$ is consistent with
$\{x\}=\alpha$ for some $\alpha\in\p{\vec a}_\Q$, we substitute $x+\alpha$ for~$x$, redo all the transformations above,
and replace $\{x\}$ everywhere with~$0$; otherwise, we can replace each equality $n\{x\}+\sum_in_ia_i=0$, $n\ne0$, with
$\bot$. This leaves only equations $\sum_in_ia_i=0$, which can be also replaced with~$\bot$ (except when
$\vec n=\vec0$) due to linear independence.

Thus, we may assume $\Gamma=\{\fii_t:t\in\N\}$ where $t\mapsto\fii_t$ is recursive, and each $\fii_t$ is a disjunction
of formulas of the form
\begin{align}
\label{eq:13}n\{x\}&>\sum_in_ia_i,\\
\label{eq:14}n\fl x&\ge\sum_in_ia_i,\\
\label{eq:15}\fl x&\equiv k\pmod m,\\
\label{eq:16}a_i&\equiv k\pmod m,
\end{align}
where $n,n_i,k,m\in\Z$, $0\le k<m$, and in \eqref{eq:16}, $a_i\in\MZ^\sM$. We will further assume that all disjuncts
\eqref{eq:15} and~\eqref{eq:16} in~$\fii_t$ use the same modulus $m=m_t$ such that $m_t\mid m_s$ whenever $t<s$: this
can be achieved by defining $m_t$ as the least common multiple of all moduli used in $\fii_0,\dots,\fii_t$, and
replacing each congruence modulo $m\mid m_t$ by an appropriate disjunction of congruences modulo~$m_t$. Moreover, we
make sure $m_t\le t$ by redefining $\fii_t$ as $\fii_{t'}$, where $t'\le t$ is maximal such that $m_{t'}\le t$.

Finally, we make sure $t\mapsto\fii_t$ (with $t$ given in unary, and the coefficients in \eqref{eq:13}--\eqref{eq:16}
written either way) is computable by a $\tc$~function. Since the function as given so far is recursive, we can write
$s=\fii_t\iff\exists r\,P(t,s,r)$, where $P$ is computable in linear time when $t,r$ are given in binary; we assume $s$
is naturally given as a binary string, which we also interpret as a G\"odel number written in binary. Then given $t$ in
unary, we can compute in $\tc$ the largest $t'\le t$ such that $\forall t''\le t'\,\exists s,r\le t\,P(t'',s,r)$ (here
we work with $t'',s,r$ as unary numbers; they have length $O(\log n)$ when converted to binary, hence $P(t'',s,r)$ can
be evaluated in logarithmic time, and therefore in~$\tc$). Taking the $s,r\le t$ such that $P(t',s,r)$ and converting
$s$ to binary, we obtain the representation of $\fii_{t'}$, which we define to be $\fii'_t$. Thus, $t\mapsto\fii'_t$ is
$\tc$-computable, and since $t\mapsto t'$ is an unbounded nondecreasing function, $\{\fii'_t:t\in\N\}$ still has all
the properties we required from $\{\fii_t:t\in\N\}$ above. Thus, we may simply assume that $t\mapsto\fii_t$ is
$\tc$-computable.

Let $T(\fii,x,\vec a)$ be a $\tc$ truth predicate for disjunctions of formulas of the form
\eqref{eq:13}--\eqref{eq:16} (with binary rational inputs $x,\vec a$), defined in the obvious way: we evaluate in
parallel each disjunct using addition, multiplication, and division with remainder (to determine $\{x\}$, $\fl x$, and
the congruences). We only need that $\sM\model\fii(x,\vec a)\eq T(\fii,x,\vec a)$ for standard formulas~$\fii$.

The final step is to construct a $\tc$~function $S(t,\vec a)$ (with $t$ in unary) that computes a rational~$x$
satisfying $\ET_{s\le t}\fii_s(x,\vec a)$, provided one exists; again, we need it to work for standard~$t$:
\begin{equation}\label{eq:17}
\forall t\in\N\:\sM\model\ET_{s\le t}\fii_s\bigl(S(t,\vec a),\vec a\bigr).
\end{equation}
Let us first observe that this will finish the proof: using~\eqref{eq:17}, the $\tc$~formula
\[\forall s\le t\,T\bigl(\fii_s,S(t,\vec a),\vec a\bigr)\]
holds in~$\sM$ for all $t\in\N$, hence it also holds for some nonstandard unary~$t$ by overspill (which follows from $\tc$-induction).
Then $u=S(t,\vec a)$ satisfies $\fii_s(u,\vec a)$ for all standard~$s$, i.e., it realizes $\Gamma(x,\vec a)$.

We define $S(t,\vec a)$ so that it works as follows:
\begin{itemize}
\item Let $V$ be the set of all rationals of the form $\frac1n\sum_in_ia_i$ such that $n\ne0$ and \eqref{eq:13} occurs
in $\fii_s$ for some $s\le t$. Sort $(V\cap(0,1))\cup\{0,1\}$ as $\{v_i:i\le p\}$, $0=v_0<v_1<\dots<v_p=1$.
\item Let $W$ consist of each $\CL{\frac1n\sum_in_ia_i}$ (for $n>0$) or $\Fl{\frac1n\sum_in_ia_i}+1$ (for $n<0$) such
that \eqref{eq:14} occurs in $\fii_s$ for some $s\le t$. Sort $W\cup\{-\infty,+\infty\}$ as $\{w_j:j\le q\}$,
$-\infty=w_0<w_1<\dots<w_q=+\infty$.
\item Determine $m_t\le t$.
\item Let $X=\bigl\{v'_i+w_{j,k}:i<p,j<q,k<m_t\bigr\}$, where $w_{j,k}=m_t\cl{m_t^{-1}(w_j-k)}+k$ for $j>0$,
$w_{0,k}=w_{1,k}-k$, and $v'_i=\frac12(v_i+v_{i+1})$.
\item Output $\min\bigl\{x\in X:\forall s\le t\,T(\fii_s,x,\vec a)\bigr\}$, if this set is nonempty.
\end{itemize}
Notice that $w_{j,k}$ is the least integer $x\ge w_j$ such that $x\equiv k\pmod{m_t}$ (for $j>0$), thus if
$[w_j,w_{j+1})$ contains an $x\equiv k\pmod{m_t}$, then $w_{j,k}$ is one such~$x$.

In order to show~\eqref{eq:17}, fix $u\in\MQ^\sM$ such that $\sM\model\ET_{s\le t}\fii_s(u,\vec a)$ and
$\{u\}\notin\p{\vec a}_\Q$. Let $i<p$, $j<q$, and $k<m_t$ be such that $\{u\}\in(v_i,v_{i+1})$,
$\fl u\in[w_j,w_{j+1})$, and $\fl u\equiv k\pmod{m_t}$. Then $u$ and $v'_i+w_{j,k}\in X$ satisfy the same formulas of
the form \eqref{eq:13}--\eqref{eq:16} that occur in $\fii_s$, $s\le t$, hence $\sM\model\ET_{s\le
t}\fii_s(v'_i+w_{j,k},\vec a)$. It follows that the set on the last line of the definition of $S(t,\vec a)$ is
nonempty, hence $S(t,\vec a)$ outputs one of its elements, which satisfies $\ET_{s\le t}\fii_s(x,\vec a)$.
\end{Pf}
\begin{Rem}\label{rem:rec-sat}
With some effort, we could generalize Theorem~\ref{thm:vtc-rec-sat} to the statement that if $\MZ^\sM\sset G\sset\MR^\sM$ is
a divisible dense subgroup, then $\p{G,\MZ^\sM,G_\ML,+,0,1,{<}}$ is recursively saturated, where $G_\ML=G\cap\MRL^\sM$.
Write $\vec a\simeq\vec b$ if $a_i\mapsto b_i$ extends to an isomorphism of ordered groups $\p{\vec a}_\Q$
and~$\p{\vec b}_\Q$. Using arguments along the lines of Theorem~3.4 and Proposition~4.1 in D'Aquino, Knight, and
Starchenko~\cite{dks}, one can show that (1) under these assumptions,
$\forall\vec a\in G\,\forall\vec b\in\MR^\sM\,\exists\vec c\in G\,\vec a,\vec b\simeq\vec a,\vec c$, and (2)
$\p{\MR^\sM,+,{<}}$ is $\omega$-homogeneous. Then given a recursive type $\Gamma(x,\vec a^\MZ,\vec a^G)$ with
$\vec a^\MZ\in\MZ^\sM$ and $\vec a^G\in(0,1)\cap G$, we find $\vec a^\MQ\in(0,1)\cap\MQ^\sM$ such that
$1,\vec a^\MQ\simeq1,\vec a^G$ by applying (1) with $\MQ^\sM$ in place of~$G$, which ensures
$\p{G,\MZ^\sM,G_\ML,+,{<},\vec a^\MZ,\vec a^G}\equiv\p{\MQ^\sM,\MZ^\sM,\MQL^\sM,+,{<},\vec a^\MZ,\vec a^\MQ}$, thus
$\Gamma(x,\vec a^\MZ,\vec a^\MQ)$ is finitely satisfiable. Fixing its realizer $c^\MQ\in\MQ^\sM$, we use (2) and~(1) to
find $\gamma^\MR\in[0,1)$ and $\gamma^G\in[0,1)\cap G$ such that
$1,\vec a^\MQ,\{c\}\simeq1,\vec a^G,\gamma^\MR\simeq1,\vec a^G,\gamma^G$, thus
$\p{G,\MZ^\sM,G_\ML,+,{<},\vec a^\MZ,\vec a^G,\fl c+\gamma^G}\equiv\p{\MQ^\sM,\MZ^\sM,\MQL^\sM,+,{<},\vec a^\MZ,\vec
a^\MQ,c}$, therefore $\fl c+\gamma^G$ realizes $\Gamma(x,\vec a^\MZ,\vec a^G)$. We leave the details to the interested
reader.

In any case, this generalization is not needed to prove the following consequence:
\end{Rem}
\begin{Cor}\label{cor:ctbl-iso-rec-sat}
If $\sM\model\vtc$ is countable, then $\p{\MQ^\sM,\MZ^\sM,+,0,1,{<}}\simeq\p{\MQL^\sM,\MZL^\sM,+,0,1,{<}}$.
Consequently, $\p{\MR^\sM,\MZ^\sM,+,0,1,{<}}\simeq\p{\MRL^\sM,\MZL^\sM,+,0,1,{<}}$, and $\MR^\sM$ expands to an
exponential field with exponential IP $\MZ^\sM$.
\end{Cor}
\begin{Pf}
We may assume $\sM$ is nonstandard. Then $\p{\MQ^\sM,\MZ^\sM,+,0,1,{<}}$ and $\p{\MQL^\sM,\MZL^\sM,+,0,1,{<}}$ are
elementarily equivalent (being models of the complete theory~$\GG$), and definable in the countable recursively
saturated structure $\p{\MQ^\sM,\MZ^\sM,\MQL^\sM,+,0,1,{<}}$, hence they are isomorphic by Theorem~\ref{thm:rec-sat}. The
restriction of any such isomorphism to~$\MZ^\sM$ extends to an isomorphism
$f\colon\p{\MR^\sM,\MZ^\sM,+,0,1,{<}}\to\p{\MRL^\sM,\MZL^\sM,+,0,1,{<}}$ by Lemma~\ref{lem:z-r}. (Alternatively, the
original isomorphism of the ordered groups $\MQ^\sM$ and~$\MQL^\sM$ extends to an isomorphism of their completions
$\MR^\sM$ and~$\MRL^\sM$, respectively, and it continues to preserve $\MZ^\sM$.) Using Theorem~\ref{thm:vtcanal},
$\exp(x)=2^{f(x)}$ defines an exponential on $\MR^\sM$ such that $\exp[\MN^\sM]\sset\MN^\sM$.
\end{Pf}

Note that Corollary~\ref{cor:ctbl-iso-rec-sat} does not yet make $\MR^\sM$ into a real-closed exponential field (even though
it is real-closed and an exponential field), as the constructed exponential need not satisfy the growth axiom
$\exp(x)>x$. We do not know how to obtain this condition using an abstract result such as Theorem~\ref{thm:rec-sat}, but as
we are going to see, it can be arranged by an adaptation of the usual back-and-forth proof of Theorem~\ref{thm:rec-sat}.
\begin{Thm}\label{thm:rcef}
If $\sM\model\vtc$ is countable, there exists an isomorphism $f\colon\p{\MZ^\sM,+,0,1,{<}}\to\p{\MZL^\sM,+,0,1,{<}}$
such that $2^{f(x)}>x$ for all $x\in\MN^\sM$. Consequently, $\MR^\sM$ expands to a real-closed exponential field with exponential IP $\MZ^\sM$.
\end{Thm}
\begin{Pf}
It suffices to prove the first part: then $\ob f(x)=f\bigl(\fl x\bigr)+\{x\}$ gives an isomorphism $\ob
f\colon\p{\MR^\sM,\MZ^\sM,+,0,1,{<}}\to\p{\MRL^\sM,\MZL^\sM,+,0,1,{<}}$ by Lemma~\ref{lem:z-r}, and it satisfies $2^{\ob
f(x)}\ge2^{f(\fl x)}\ge\fl x+1>x$ for $x\ge0$, hence we can apply Corollary~\ref{cor:r-rl-eip}.

Let ${\log}\colon\MR_{>0}^\sM\to\MRL^\sM$ denote the inverse of $2^x$, and $\lang$ the language of ordered groups; if
$\vec a\in\MZ^\sM$ and $\vec b\in\MZL^\sM$ have the same length, we write $\vec a\equiv_\lang\vec b$ for
$\p{\MZ^\sM,+,{<},\vec a}\equiv\p{\MZL^\sM,+,{<},\vec b}$. Fix enumerations $\MZ^\sM=\{u_n:n\in\N\}$ and
$\MZL^\sM=\{v_n:n\in\N\}$. By induction on~$n$, we will define sequences $\{a_n:n\in\N\}\sset\MZ^\sM$ and
$\{b_n:n\in\N\}\sset\MZL^\sM$ with the following properties:
\begin{enumerate}
\item\label{item:8} $a_0=b_0=1$, $a_{2n+1}=u_n$, and $b_{2n+2}=v_n$.
\item\label{item:9} $\vec a_{<n}\equiv_\lang\vec b_{<n}$, where $\vec a_{<n}=\p{a_i:i<n}$, and similarly for~$\vec
b_{<n}$.
\item\label{item:10} For all $\vec q\in\Q^n$,  $\sum_{i<n}q_ia_i>0\implies\sum_{i<n}q_ib_i>\log\sum_{i<n}q_ia_i$.
\end{enumerate}
Notice that by Presburger quantifier elimination, \ref{item:9} is (in view of $a_0=b_0=1$) equivalent to
\begin{enumerate}
\item[(ii$'$)] $\sum_{i<n}q_ia_i\leqg0\iff\sum_{i<n}q_ib_i\leqg0$ for all $\vec q\in\Q^n$ and
${\leqg}\in\{{<},{=},{>}\}$, and $a_i\equiv b_i\pmod m$ for all $i<n$ and $m\in\N_{>0}$.
\end{enumerate}
Moreover, \ref{item:10} is equivalent to
\begin{enumerate}
\item[(iii$'$)] For all $\vec q\in\Q^n$, $\sum_{i<n}q_ia_i>\N\implies\sum_{i<n}q_ib_i>\N\log\sum_{i<n}q_ia_i$\,:
\end{enumerate}
since $\sum_iq_ia_i\in\N^{-1}\MZ^\sM$, either $\sum_iq_ia_i\in\Q$ is standard, in which case (using~\ref{item:9})
$\sum_iq_ib_i=\sum_iq_ia_i>\log\sum_iq_ia_i$ holds automatically, or $\sum_iq_ia_i>\N$, in which
case \ref{item:10} implies
\[\sum_{i<n}q_ib_i=2k\sum_{i<n}\frac{q_i}{2k}b_i>2k\log\sum_{i<n}\frac{q_i}{2k}a_i
  =2k\Bigl(\log\sum_{i<n}q_ia_i-\log 2k\Bigr)>k\log\sum_{i<n}q_ia_i\]
for all $k\in\N_{>0}$.

As indicated by~\ref{item:8}, we put $a_0=b_0=1$, which satisfies \ref{item:8}--\ref{item:10} by the discussion above.
Assume that $n>0$, and $\vec a_{<n}$ and $\vec b_{<n}$ have been defined such that \ref{item:8}--\ref{item:10} hold; we
will define $a_n$ and~$b_n$.

If $n$ is odd, we put $a_n=u_{(n-1)/2}$. If $a_n=\sum_{i<n}q_ia_i$ for some $\vec q\in\Q^n$, we define
$b_n=\sum_{i<n}q_ib_i$; then \ref{item:8}--\ref{item:10} follow from the induction hypothesis. If $a_n\notin\p{\vec
a_{<n}}_\Q$, we use Lemma~\ref{lem:fg-cof} to find $a\in\p{\vec a_{\le n}}_\Q$ such that $\N^{-1}a$ is cofinal above
$\bigl\{z:\forall x\in\p{\vec a_{\le n}}_\Q\,\bigl(\abs x\le\abs z\to x\in\p{\vec a_{<n}}_\Q\bigr)\bigr\}$, that is,
\begin{equation}\label{eq:18}
\forall x\in\p{\vec a_{\le n}}_\Q\bez\p{\vec a_{<n}}_\Q\:\exists k\in\N\:k\abs x\ge a.
\end{equation}
We may assume $a\in\MZ^\sM$. Observe $a>\N$. Putting $c=\fl{\log a}$, we claim that the $\lang_\GGG$-type
\begin{align*}
\Gamma(x)&=\Bigl\{x\lessgtr\sum_{i<n}q_ib_i\eq a\lessgtr\sum_{i<n}q_ia_i:\vec q\in\Q^n,{\lessgtr}\in\{{<},{>}\}\Bigr\}\\
&\qquad{}\cup\bigl\{x\equiv a\pmod m:m\in\N_{>0}\bigr\}\cup\{x>kc:k\in\N\}\cup\{L(x)\}\\
&\equiv\Bigl\{x\lessgtr\sum_{i<n}q_ib_i:\vec q\in\Q^n,{\lessgtr}\in\{{<},{>}\},a\lessgtr\sum_{i<n}q_ia_i\Bigr\}\\
&\qquad{}\cup\bigl\{x\equiv a\pmod m:m\in\N_{>0}\bigr\}\cup\{x>k\log a:k\in\N\}\cup\{x\in\MZL^\sM\}
\end{align*}
is finitely satisfiable. If a finite $\Gamma'\sset\Gamma$ involves no linear inequality $x<\sum_iq_ib_i$, it is
satisfied by any sufficiently large $x\in\MZL^\sM$ satisfying the congruences. Otherwise, it is equivalent to
\begin{equation}\label{eq:19}
\max\Bigl\{\sum_{i<n}r_ib_i,k\log a\Bigr\}<x<\sum_{i<n}q_ib_i\land x\equiv a\pmod m
\end{equation}
for some $\vec q,\vec r\in\Q^n$ and $k,m\in\N_{>0}$ such that $\sum_ir_ia_i<a<\sum_iq_ia_i$. Notice that
\[\sum_{i<n}q_ia_i-\sum_{i<n}r_ia_i>\sum_{i<n}q_ia_i-a>\N\colon\]
if not, then using $\p{\vec a_{\le n}}_\Q\sset\N^{-1}\MZ^\sM$ we obtain that $\sum_iq_ia_i-a\in\Q$, hence
$a\in\p{\vec a_{<n},1}_\Q=\p{\vec a_{<n}}_\Q$, a contradiction. Thus, $\sum_iq_ib_i-\sum_ir_ib_i>\N$ as well, using
\ref{item:9} of the induction hypothesis. Likewise, the induction hypothesis gives $\sum_iq_ib_i>\N\log\sum_iq_ia_i$,
thus $\sum_iq_ib_i>k\log a+\N$. It follows that the interval defined by the bounds in~\eqref{eq:19} has nonstandard
length, and as such contains an element satisfying the congruence.

Using Theorem~\ref{thm:vtc-rec-sat}, $\Gamma(x)$ is realized by an element~$b\in\MZL^\sM$. Clearly, $\vec
a_{<n},a\equiv_\lang\vec b_{<n},b$. We claim that
\begin{equation}\label{eq:20}
\alpha:=\sum_{i<n}q_ia_i+qa>\N\implies\beta:=\sum_{i<n}q_ib_i+qb>\N\log\alpha
\end{equation}
for all $\vec q\in\Q^n$, $q\in\Q$. If $q=0$, this follows from the induction hypothesis. If $\alpha>\N a$, we have
$\frac12\sum_iq_ia_i<\alpha<2\sum_iq_ia_i$, and $2\sum_iq_ib_i>\N\log\bigl(2\sum_iq_ia_i\bigr)$ by the induction
hypothesis, thus $\beta>\frac12\sum_iq_ib_i>\N\log\alpha$. In the remaining case, \eqref{eq:18} shows that
$\frac1ka<\alpha<ka$ for some $k\in\N$, thus $\beta>\frac1kb>\N\log a$ by the definition of~$\Gamma$, and
$\beta>\N\log\alpha$.

Since $a\in\p{\vec a_{\le n}}_\Q\bez\p{\vec a_{<n}}_\Q$, we have $a_n=\sum_{i<n}q_ia_i+qa$ for some $\vec q\in\Q^n$,
$q\in\Q$. Putting $b_n=\sum_iq_ib_i+qb$, the condition $\vec a_{<n},a\equiv_\lang\vec b_{<n},b$ implies that
$b_n\in\MZL^\sM$ and $\vec a_{\le n}\equiv_\lang\vec b_{\le n}$, and \eqref{eq:20} implies \ref{item:10}.

Now, let $n$ be even, and put $b_n=v_{n/2-1}$. As in the previous case, we may assume $b_n\notin\p{\vec b_{<n}}_\Q$,
and we can fix $b\in\MZ_{\ML,>0}^\sM\cap\p{\vec b_{\le n}}_\Q\bez\p{\vec b_{<n}}_\Q$ such that
\[\forall x\in\p{\vec b_{\le n}}_\Q\bez\p{\vec b_{<n}}_\Q\:\exists k\in\N\:k\abs x\ge b.\]
The same argument as above shows that any realizer $a\in\MZ^\sM$ of the type
\begin{align*}
\Gamma(x)&=\Bigl\{x\lessgtr\sum_{i<n}q_ia_i\eq b\lessgtr\sum_{i<n}q_ib_i:\vec q\in\Q^n,{\lessgtr}\in\{{<},{>}\}\Bigr\}\\
&\qquad{}\cup\bigl\{x\equiv b\pmod m:m\in\N_{>0}\bigr\}\cup\{x<2^{b/k}:k\in\N\}
\end{align*}
can be used to construct $a_n\in\MZ^\sM$ such that \ref{item:9} and~\ref{item:10} hold. The catch is that there is no
obvious way how to define $\{x<2^{b/k}:k\in\N\}$ using only finitely many parameters, hence we need to replace this
part.

Using Lemma~\ref{lem:fg-cof}, there is $d=\sum_iq_ib_i\in\MZ_{\ML,>0}^\sM$, $\vec q\in\Q^n$, such that $\N d$ is cofinal
in the convex subgroup $C=\{x\in\p{\vec b_{<n}}_\Q:\exists k\in\N\,\abs x\le kb\}$ of $\p{\vec b_{<n}}_\Q$. Put
$c=\sum_iq_ia_i$. Since $b\ge\frac1kd$ for some $k\in\N$, and $d>\N\log c$ by the induction hypothesis, we have
$\N^{-1}b>2\log c$, thus $\Gamma(x)$ is implied by the type
\begin{align*}
\Gamma'(x)&=\Bigl\{x\lessgtr\sum_{i<n}q_ia_i\eq b\lessgtr\sum_{i<n}q_ib_i:\vec q\in\Q^n,{\lessgtr}\in\{{<},{>}\}\Bigr\}\\
&\qquad{}\cup\bigl\{x\equiv b\pmod m:m\in\N_{>0}\bigr\}\cup\{x<c^2\}.
\end{align*}
It remains to verify that $\Gamma'(x)$ is finitely satisfiable. As before, this amounts to showing that if
$\sum_iq_ib_i<b<\sum_ir_ib_i$ and $m\in\N_{>0}$, there exists $x\in\MZ^\sM$ such that
\[\sum_{i<n}q_ia_i<x<\min\Bigl\{\sum_{i<n}r_ia_i,c^2\Bigr\}\land x\equiv b\pmod m,\]
which in turn holds if the difference between the two bounds is nonstandard. Also, $b\notin\p{\vec b_{<n}}_\Q$ again
implies that $\sum_ir_ib_i-\sum_iq_ib_i>\N$ and $\sum_ir_ia_i-\sum_iq_ia_i>\N$. Finally, the cofinality of $\N d$
in~$C$ ensures that $\sum_iq_ib_i\le kd$ for some $k\in\N$, thus $\sum_iq_ia_i\le kc<c^2$.
\end{Pf}

As we indicated in Section~\ref{sec:real-expon-models}, the proof of Theorem~\ref{thm:rcef} only used the recursive
saturation of $\p{\MZ^\sM,\MZL^\sM,+,0,1,{<}}$. Moreover, we used the predicate $L(x)$ only in a very limited way,
namely to realize a type of the form $\Gamma(x)\cup\{L(x)\}$ where $\Gamma$ is in the language of ordered groups. It
would not be difficult to eliminate it entirely, so that the proof would only use the recursive saturation of the
Presburger reduct $\p{\MZ^\sM,+,0,1,{<}}$. While this would not significantly simplify the proof of
Theorem~\ref{thm:vtc-rec-sat} proper, we could dispense with the material in Sections \ref{sec:theory-three-groups}
and~\ref{sec:recurs-satur-models} in favour of the standard quantifier elimination for $\Z$-groups. However, as
we already stressed, we consider the full statement of Theorem~\ref{thm:vtc-rec-sat} to be intrinsically interesting in
its own right, and therefore keep all the results.

The proof of Theorem~\ref{thm:rcef} essentially relies on the countability of~$\sM$ (unlike Corollary~\ref{cor:idz+exp}, which
applies to arbitrarily large models). We do not know to what extent it can be generalized to uncountable models, though
we can at least infer the following on general principle:
\begin{Cor}\label{cor:rcef}
Every model $\sM_1\model\vtc$ has an elementary extension $\sM$ of the same cardinality that satisfies the conclusions
of Theorem~\ref{thm:rcef}.
\end{Cor}
\begin{Pf}
Let $\sM_0$ be a countable elementary submodel of~$\sM_1$. By Theorem~\ref{thm:rcef}, there exists an isomorphism
$f_0\colon\p{\MZ^{\sM_0},+,0,1,{<}}\to\p{\MZL^{\sM_0},+,0,1,{<}}$ such that $2^{f_0(x)}>x$ for all $x\in\MN^{\sM_0}$.
Since $\Th(\sM_0,f_0)$ is consistent with the elementary diagram of~$\sM_1$, there exists an elementary extension
$\sM$ of~$\sM_1$ of the same cardinality and a function~$f$ such that $\p{\sM,f}\equiv\p{\sM_0,f_0}$, which ensures
that $f$ an isomorphism $\p{\MZ^\sM,+,0,1,{<}}\to\p{\MZL^\sM,+,0,1,{<}}$ satisfying $2^{f(x)}>x$ for all
$x\in\MN^\sM$. Then $\MR^\sM$ expands to a real-closed exponential field with exponential IP $\MZ^\sM$ by
Lemma~\ref{lem:z-r} and Corollary~\ref{cor:r-rl-eip}.
\end{Pf}

\section{Conclusion and open problems}\label{sec:conclusion}

We have shown that countable models of~$\vtc$ are exponential IP of real-closed exponential fields; among other things,
this severely limits the first-order consequences of being an exponential IP of a RCEF. Our work suggests various
follow-up problems. The first one is that we could not prove much of anything about uncountable models of~$\vtc$,
besides the rather unsatisfactory Corollary~\ref{cor:rcef}:
\begin{Que}\label{que:uncount}
Is every uncountable model of $\vtc$ an exponential IP of a real-closed exponential field? If not, can we characterize
the models that are?
\end{Que}
Let us also recall a question from~\cite{ej-kol:realclosures}: are real-closed fields with IP $\sM\model\vtc$, $\sM$
nonstandard, recursively saturated?

We may also look at other theories. Due to the $\tc$-completeness of integer multiplication, $\vtc$ is the weakest
reasonable theory in the setup of Zambella-style two-sorted theories of arithmetic whose models carry a ring structure,
as $\vtc$ is axiomatizable by the totality of multiplication over the standard base theory~$\thry{V^0}$. But of course,
we may consider weaker or incomparable theories in the basic one-sorted language of arithmetic.

In particular, additive reducts of nonstandard models of~$\thry{IE_1}$ are recursively saturated due to
Wilmers~\cite{wilmers}, and the corresponding property of~$\vtc$ was one of the main ingredients of our proof of
Theorem~\ref{thm:rcef}. On the other hand, it is unclear if we can complement this for every $\sM\model\thry{IE_1}$ with a
construction of an exponential $2^x\colon\p{L,+,{<}}\simeq\p{\MR^\sM_{>0},\cdot,{<}}$ for a convex subgroup
$L\sset\MR^\sM$ (satisfying \ref{item:7} of Theorem~\ref{thm:3g-rec-sat}), which was the other main ingredient. Notice that
$\thry{IE_1}$ (or even $\idz$) has nonstandard models~$\sM$ that are polynomially bounded in the sense that
$\{a^n:n\in\N\}$ is cofinal in~$\sM$ for some element~$a$; it is easy to see that such models cannot be exponential IP
of any exponential fields.
\begin{Que}\label{que:ie1}
Is every non-polynomially-bounded countable model of~$\thry{IE_1}$ (or at least, $\idz$) an exponential IP of a RCEF?
\end{Que}

Concerning first-order consequences of being an exponential IP of a RCEF, the author is not actually aware of any
whatsoever beside the obvious ones, which suggests:
\begin{Que}\label{que:iopen}
Does every model of $\io$ have an elementary extension to an exponential IP of a RCEF?
\end{Que}
We observe that every model of~$\io$ has an elementary extension that is a (not necessarily exponential) IP of a RCEF
by a simple application of Robinson's joint consistency theorem.

\section*{Acknowledgement}

I want to thank the anonymous referee for useful suggestions.

The research was supported by grant 23-04825S of GA \v CR. The Institute of Mathematics of the Czech Academy of
Sciences is supported by RVO: 67985840.


\begin{thebibliography}{10}
\bibliographyhook

\bibitem{founif}
David A.~Mix Barrington, Neil Immerman, and Howard Straubing, \emph{On
  uniformity within {$\mathit{NC}^1$}}, Journal of Computer and System Sciences
  41 (1990), no.~3, pp.~274--306.

\bibitem{bar-sch}
Jon Barwise and John Schlipf, \emph{An introduction to recursively saturated
  and resplendent models}, Journal of Symbolic Logic 41 (1976), no.~2,
  pp.~531--536.

\bibitem{buss}
Samuel~R. Buss, \emph{Bounded arithmetic}, Bibliopolis, Naples, 1986, revision
  of 1985 Princeton University Ph.D. thesis.

\bibitem{cdk}
Merlin Carl, Paola D'Aquino, and Salma Kuhlmann, \emph{On the value group of a
  model of {P}eano arithmetic}, Forum Mathematicum 29 (2017), no.~4,
  pp.~951--957.

\bibitem{carl-krapp}
Merlin Carl and Lothar~Sebastian Krapp, \emph{Models of true arithmetic are
  integer parts of models of real exponentiation}, Journal of Logic and
  Analysis 13 (2021), no.~3, pp.~1--21.

\bibitem{cmw}
Patrick C{\'e}gielski, Kenneth McAloon, and George Wilmers, \emph{Mod{\`e}les
  r{\'e}cursivement satur{\'e}s de l'addition et de la multiplication des
  entiers naturels}, in: Logic {C}olloquium '80 (D.~van Dalen et~al., eds.),
  Studies in Logic and the Foundations of Mathematics vol. 108, North-Holland,
  1982, pp.~57--68.

\bibitem{cook-ngu}
Stephen~A. Cook and Phuong Nguyen, \emph{Logical foundations of proof
  complexity}, Perspectives in Logic, Cambridge University Press, New York,
  2010.

\bibitem{dkkl:rcef}
Paola D'Aquino, Julia~F. Knight, Salma Kuhlmann, and Karen Lange, \emph{Real
  closed exponential fields}, Fundamenta Mathematicae 219 (2012), pp.~163--190.

\bibitem{dks}
Paola D'Aquino, Julia~F. Knight, and Sergei Starchenko, \emph{Real closed
  fields and models of {P}eano arithmetic}, Journal of Symbolic Logic 75
  (2010), no.~1, pp.~1--11.

\bibitem{tc0}
Andr{\'a}s Hajnal, Wolfgang Maass, Pavel Pudl{\'a}k, M{\'a}ri{\'o} Szegedy, and
  Gy{\"o}rgy Tur{\'a}n, \emph{Threshold circuits of bounded depth}, Journal of
  Computer and System Sciences 46 (1993), no.~2, pp.~129--154.

\bibitem{jen-ehr:prob}
Don Jensen and Andrzej Ehrenfeucht, \emph{Some problem in elementary
  arithmetics}, Fundamenta Mathematicae 92 (1976), no.~3, pp.~223--245.

\bibitem{ej:vtc0iopen}
Emil Je{\v r}{\'a}bek, \emph{Open induction in a bounded arithmetic for
  {$\mathrm{TC}^0$}}, Archive for Mathematical Logic 54 (2015{\gobble a}),
  no.~3--4, pp.~359--394.

\bibitem{ej:vtcimul}
\bysame, \emph{Iterated multiplication in {$\mathit{VTC}^0$}}, Archive for
  Mathematical Logic 61 (2022{\gobble b}), no.~5--6, pp.~705--767.

\bibitem{ej:vtcanal}
\bysame, \emph{Elementary analytic functions in {$\mathsf{VTC^0}$}}, Annals of
  Pure and Applied Logic 174 (2023{\gobble b}), no.~6, article no.~103269,
  50~pp.

\bibitem{ej-kol:realclosures}
Emil Je{\v r}{\'a}bek and Leszek~A. Ko{\l}odziejczyk, \emph{Real closures of
  models of weak arithmetic}, Archive for Mathematical Logic 52 (2013{\gobble
  a}), no.~1--2, pp.~143--157.

\bibitem{joh-pol:d1cr}
Jan Johannsen and Chris Pollett, \emph{On the {$\Delta^b_1$}-bit-comprehension
  rule}, in: Logic {C}olloquium '98: Proceedings of the 1998 {ASL} {E}uropean
  {S}ummer {M}eeting held in {P}rague, {C}zech {R}epublic (S.~R. Buss,
  P.~H{\'a}jek, and P.~Pudl{\'a}k, eds.), ASL, 2000, pp.~262--280.

\bibitem{krapp:phd}
Lothar~Sebastian Krapp, \emph{Algebraic and model theoretic properties of
  o-minimal exponential fields}, Ph.D. thesis, Universit{\"a}t Konstanz, 2019.

\bibitem{lessan:phd}
Hamid Lessan, \emph{Models of arithmetic}, Ph.D. thesis, University of
  Manchester, 1978.

\bibitem{mou-res}
Marie-H{\'e}l{\`e}ne Mourgues and Jean-Pierre Ressayre, \emph{Every real closed
  field has an integer part}, Journal of Symbolic Logic 58 (1993), no.~2,
  pp.~641--647.

\bibitem{ngu-cook}
Phuong Nguyen and Stephen~A. Cook, \emph{Theories for {$\mathit{TC}^0$} and
  other small complexity classes}, Logical Methods in Computer Science 2
  (2006), no.~1, article no.~3, 39~pp.

\bibitem{par-sch}
Ian Parberry and Georg Schnitger, \emph{Parallel computation with threshold
  functions}, Journal of Computer and System Sciences 36 (1988), no.~3,
  pp.~278--302.

\bibitem{ress:eip}
Jean-Pierre Ressayre, \emph{Integer parts of real closed exponential fields},
  in: Arithmetic, proof theory, and computational complexity (P.~Clote and
  J.~Kraj{\'\i}{\v c}ek, eds.), Oxford Logic Guides vol.~23, Oxford University
  Press, 1993, pp.~278--288.

\bibitem{scott-cof}
Dana Scott, \emph{On completing ordered fields}, in: Applications of Model
  Theory to Algebra, Analysis, and Probability (W.~A.~J. Luxemburg, ed.), Holt,
  Rinehart and Winston, New York, 1969, pp.~274--278.

\bibitem{sheph}
John~C. Shepherdson, \emph{A nonstandard model for a free variable fragment of
  number theory}, Bul\-le\-tin de l'Aca\-d{\'e}\-mie Po\-lo\-naise des
  \hbox{Sciences}, S{\'e}\-rie des \hbox{Sciences} Ma\-th{\'e}\-ma\-tiques,
  As\-tro\-no\-miques et Phy\-siques 12 (1964), no.~2, pp.~79--86.

\bibitem{wer:topf}
Seth Warner, \emph{Topological fields}, North-Holland Mathematics Studies vol.
  157, North-Holland, New York, 1989.

\bibitem{wilmers}
George Wilmers, \emph{Bounded existential induction}, Journal of Symbolic Logic
  50 (1985), no.~1, pp.~72--90.

\bibitem{zamb:notes}
Domenico Zambella, \emph{Notes on polynomially bounded arithmetic}, Journal of
  Symbolic Logic 61 (1996), no.~3, pp.~942--966.

\end{thebibliography}
 \providecommand\gobble[1]{} \ifx\url\undefined {\catcode`\/=13
  \gdef/{\string/\futurelet\nexttoken\finishslash}
  \gdef\finishslash{\ifx\nexttoken/\else\penalty\relpenalty\fi}}
  \def\url{\begingroup\catcode`\~=12 \catcode`\_=12 \catcode`\/=13 \finishurl}
  \def\finishurl#1{\texttt{#1}\endgroup} \fi \providecommand\href[2]{\url{#2}}
  \providecommand\dotminus{\mathbin{\scriptstyle\dot{\smash{\textstyle-}}}}
  \providecommand\hyph{\nobreak\hskip0pt-\hskip0pt\relax}
\providecommand\bysame{\leavevmode\hbox to5em{\hrulefill}\thinspace}
\providecommand\bibliographyhook{}

\end{document}